\documentclass[a4paper,12pt]{article}
\usepackage{amssymb,amsmath,amsthm,latexsym}
\usepackage[croatian,english]{babel}

\usepackage{color}

\newcommand{\normal}{\color{black}}

\usepackage[a4paper,margin=1in]{geometry} %simple margins
\usepackage{verbatim}%for the command "comment" to comment out whole parts
\usepackage{epsf,graphicx}

\theoremstyle{plain}
\newtheorem{thm}{Theorem}[section]

\newtheorem{cor}[thm]{Corollary}
\newtheorem{lemma}[thm]{Lemma}

\theoremstyle{definition}
\newtheorem{defn}[thm]{Definition}

\numberwithin{equation}{section}

\newcommand{\R}{{\mathbb R}}
\newcommand{\N}{{\mathbb N}}
\newcommand{\Z}{{\mathbb Z}}

\newcommand{\E}{\mathbb{E}\,}

\def\P{{\mathbb P}}

\begin{document}

\title{\Large\bfseries Ruin probability for discrete risk processes}

\author{\itshape
    Ivana Ge\v{c}ek Tu{\dj}en\\
    \\
    Department of Mathematics\\
    University of Zagreb, Zagreb, Croatia\\
    Email: igecek@math.hr
   \normal
}

\date{}

\maketitle

\begin{abstract}
\noindent
We study the discrete time risk process modelled by the skip-free random walk and we derive the results connected to the ruin probability, such as crossing the fixed level, for this kind of process. We use the method relying on the classical ballot theorems to derive these results and compare them to the results obtained for the continuous time version of the risk process. We further generalize this model by adding the perturbation and, still relying on the skip-free structure of that process, we generalize the previous results on crossing the fixed level for the generalized discrete time risk process.


\end{abstract}

\section{Introduction}\label{sec-1}
In the classical ruin theory one usually observes the risk process
\begin{equation}\label{1:eq.1} X(t)=ct-\sum_{i=1}^{N(t)}Y_i~,~~t\geq 0~,\end{equation}
 where $c>0$ represents the premium rate (we assume that there are the incoming premiums which arrive from the policy holders),
$(Y_i~:~i\in \N)$ is an i.i.d. sequence of nonnegative random variables with common distribution $F$ (which usually represent the policy holders' claims) and $(N(t)~:~t\geq 0)$ a homogeneous Poisson process of rate $\lambda>0$, independent of $(Y_i~:~i\in \N)$.
This basic process was generalized by many authors and we will follow the approach used in \cite{HPSV1}, \cite{HPSV2} and \cite{IGT}, which means that in the continuous  time case one observes the generalized risk process
\begin{equation}\label{1:eq.2}X(t)=ct-C(t)+Z(t)~,~~t\geq 0~,\end{equation}
where $(C(t)~:~t\geq 0)$ is a subordinator and $(Z(t)~:~t\geq 0)$ an independent spectrally negative L\'{e}vy
process. The overall process $X$ then also has the nice structure of the spectrally negative L\'evy process and the results from the fluctuation theory may be used to analyze it.
One of the main questions that is observed in this model is the question of the \emph{ruin probability}, given some initial capital
$u>0$, i.e.
\begin{equation}\label{1:eq.3}\vartheta(u)=\P(u+X(t)<0~,~~\textrm{for some $t>0$})~.\end{equation}
Furthermore, the question of the distribution of the supremum of the dual process $\widehat{X}=-X$ is of the main interest, as well as the question directly connected to it, i.e. the first passage over some fixed level. Results for the above questions can be obtained using different approaches, such as decomposing the supremum of the dual of the generalized risk process $\widehat{X}$ or Laplace transform approach, and in \cite{HPSV2} authors use the famous Tak\'acs formula in the continuous time case.\\
\\
More precisely, for $m$ independent subordinators $C_1,\ldots,C_m$ without drift and with L\'evy measures $\Lambda_1,\ldots,\Lambda_m$ such that $\E(C_i(1))<\infty$, $i=1,2,\ldots,m$, one observes the risk process
\begin{equation}\label{1:eq.4}X(t)=ct-C(t)~,~~t\geq 0~,\end{equation}
for $C=C_1+\cdots +C_m$ and $c>\E(C_1(1))+\cdots +\E(C_m(1))$ (standard net profit assumption). In \cite{HPSV2} the following result was achieved:
\begin{equation}\label{1:eq.5}\P(\widehat{\tau}_0<\infty,\widehat{X}(\widehat{\tau}_0-)\in dy,\widehat{X}(\widehat{\tau}_0)\in dx,\Delta C(\widehat{\tau}_0)=\Delta C_i(\widehat{\tau}_0))=\frac{1}{c}\Lambda_i(-y+dx)dy~,\end{equation}
for $x>0$, $y\leq 0$ and $\widehat{X}(0)=0$, $\widehat{\tau}_0=\inf\{t\geq 0~:~\widehat{X}(t)>0\}$, $\Delta C(t)=C(t)-C(t-)$, $\Delta C_i(t)=C_i(t)-C_i(t-)$, $i=1,\ldots ,m$.\\
\\
Here authors interpret processes $C_i$ as independent risk portfolios competing to cause ruin and the above formula gives the probability that the ruin will be caused by one individual portfolio. Authors further generalize the problem by adding the L\'evy process $Z$ with no positive jumps in the model (this is called the perturbed model) and achieve the similar formula.\\
\\
The focus of this paper is rather on the method which led to above result, namely, the before mentioned Tak\'acs "magic" formula. In the continuous time case, this formula can be expressed in the following way (for details see \cite{Tak}).
\begin{lemma} For the process $X$ defined as in (\ref{1:eq.4}), with $\widehat{X}(0)=0$,
\begin{equation}\label{1:eq.6}
\P\big(\sup_{0\leq s\leq t}\widehat{X}(s)>0~|~\widehat{X}(t)\big)=1-\big(-\frac{\widehat{X}(t)}{ct}\big)~.
\end{equation}
\end{lemma}
But, this result naturally arises in the discrete time case in the view of the well known ballot theorems, so the main aim of this paper is to discover how and which results for ruin probability can be obtained, using the above method, when we observe a discrete time risk process. To establish the connection with the continuous time model, we will model our discrete time risk process with the \emph{upwards skip-free} (or \emph{right continuous}) random walk, i.e. a random walk with increments less or equal than $1$. These random walks can be observed as a discrete version of the spectrally negative L\'evy processes, i.e. the processes with no positive jumps. The main connection, which is in the focus of this paper, between this discrete and the continuous model is that in both cases we are able to control the one side jumps of the process. More precisely, skip-free random walks can cross levels from one side with the jumps of any size and on the other side they can only have unit jumps.\\
\\
In this surrounding we will prove the main results of the paper and the main tool for it will be the following result (details for this type of result can be found in \cite{Tak} or \cite{Dwa}).
\begin{lemma}
Let $\xi_1,\ldots,\xi_n$ be the random variables in $\{\ldots,-3,-2,-1,0,1\}$ with cyclically interchangeable increments. Let $R(i)=\xi_1+\cdots+\xi_n$, $1\leq i\leq n$, $R(0):=0$. Then for each $0\leq k\leq n$
\begin{equation}\label{1.eq.7}
\P(R(i)>0~\textrm{for each  $1\leq i\leq n$}~|~R(n)=k)=\frac{k}{n}~.
\end{equation}
\end{lemma}
Using the skip-free structure of the random walks that model our risk process, Lemma 1.2. and some auxiliary results following from the ballot theorems (such as \emph{Kemperman's formula}, which will be explained in details in Section 2), we will derive the following main results of this paper.
\begin{thm} Let $C^1$ and $C^2$ be two independent random walks with nondecreasing increments and $\mu^i:=\E(C^i(1))<\infty$, $i=1,2$. Let $C:=C^1+C^2$, $\mu=\mu^1+\mu^2$ and
\begin{equation}\label{1:eq.8}
X(n)=n-C(n)~,~~n\geq 0~,
\end{equation}
and let us assume that $\E (X(1))>0$, i.e. $\mu<1$. Then
\begin{equation}\label{1:eq.9}
\P(\widehat{\tau_0}<\infty,\widehat{X}(\widehat{\tau_0}-1)=y,
\widehat{X}(\widehat{\tau_0})\geq x,\Delta
C^i(\widehat{\tau_0})=x+1-y)=\P(C^i(1)=x+1-y)~,
\end{equation}
for $y\leq 0$, $x>0$, $\widehat{\tau}_0=\inf\{n\geq 0~:~\widehat{X}(n)>0\}$, $\Delta C^i(n)=C^i(n)-C^i(n-1)$ (for $i=1,2$) and $\Delta C(n)=C(n)-C(n-1)$, $n\geq 0$.
\end{thm}
The above result will be generalized for $m$ independent random walks with nondecreasing increments $C^1,\ldots,C^m$, $m\in\N$, and $C=C^1+\cdots +C^m$ in the standard way.\\
\\
When we generalize the above model by adding the perturbation modelled by an \emph{upwards skip-free} random walk (or \emph{right continuous} random walk ) $Z$, i.e. the random walk with increments less or equal than $1$, we observe the perturbed discrete time risk process
\begin{equation}\label{1:eq.10}
X(n)=-C(n)+Z(n)~,~~n\geq 0~.
\end{equation}
Under the assumption that $\E(X(1))>0$ (so $\widehat{X}\to -\infty$) and with the same notation used as in the previous theorem, we will derive the following result.
\begin{thm}\label{1.eq.11}
\begin{align*}\P(\widehat{\tau}_0<\infty,&\widehat{X}(\widehat
{\tau}_0-1)=y,\widehat{X}(\widehat{\tau}_0)\geq x,\textrm{the new supremum was caused by the process $C$})\\
&=\P(C(1)\geq
x+1-y)\cdot\P(Z(1)=-1)+\P(C(1)\geq x-y)\cdot\P(Z(1)\geq
0)~.\end{align*}
\end{thm}
This result can also be generalized so that we observe the probability that the random walk $C^i$ causes the ruin ($i\in\{1,\ldots,m\}$ ), again in the standard way.

\section{Auxiliary results}\label{sec-2}

\begin{defn}
Let $(S(n)~:~n\geq 0)$ be a random walk with integer-valued increments $(Y(i)~:~i\geq 0)$, i.e. $S(n)=\sum_{i=1}^nY(i)$, $n\geq 0$, $S(0)=0$. We say that $S$ is an \emph{upwards skip-free} (or \emph{right continuous}) random walk if $\P(Y(i)\leq 1)=1$ (i.e. it's increments achieve values greater than $1$ with zero probability).\\
If $\P(Y(i)\geq -1)=1$ we say that $S$ is a \emph{downwards skip-free} (or \emph{left continuous}) random walk.
\end{defn}
To prove the results for the ruin probability for the skip-free class of random walks we will need some auxiliary results. Inspired by the approach used in \cite{HPSV2} for the continuous time case (i.e. spectrally negative L\'evy processes), we will use the results following from the famous \emph{ballot theorems}, first one dating to 1887. and formulated by \emph{Bertrand}. More precisely, let us assume that there are two candidates in the voting process in which there are $n$ voters. Candidate $A$ scores $a$ votes for the win over the candidate $B$ which scores $b$ votes, $a\geq b$. Then the probability that throughout the counting the number of votes registered for $A$ is always greater than the number of votes registered for the candidate $B$ is equal to $\frac{a-b}{a+b}=\frac{a-b}{n}$. This result was further generalized by \emph{Barbier} and proved in that generalized form by \emph{Aeplli}. Later this result was also proved by  \emph{Dvoretzky} and \emph{Motzkin} using the \emph{cyclic lemma}, which is the approach similar to the one followed in this paper.\\
\\
The main property that lies in the heart of those type of theorems is having some kind of cyclic structure.
\begin{defn}
 For random variables $\xi_1,\ldots,\xi_n$ we say that they are \emph{interchangeable} if for each
 $(r_1,\ldots ,r_n)\in \R^n$
and all permutations $\sigma$ of $\{1,2,\ldots ,n\}$,
\begin{equation}\label{2:eq.1}\P(\xi_i\leq
r_i~\textrm{for each $1\leq i\leq n$}~)=\P(\xi_i\leq
r_{\sigma(i)}~\textrm{for each $1\leq i\leq n$}~)~.
\end{equation}
For $\xi_1,\ldots
,\xi_n$ we say that they are \emph{cyclically interchangeable} if (\ref{2:eq.1}) is valid for each cyclic permutation $\sigma$ of
$\{1,2,\ldots ,n\}$.
\end{defn}
In other words, the random variables are cyclically interchangeable if their distribution law is invariant under cyclic permutations and interchangeable (in some literature, for example see \cite{Dwa}, this property is also called \emph{exchangeable}) if it is invariant under all permutations. From the definition it is clear that interchangeable variables are also cyclically interchangeable and the converse is not true.\\
\\
One version of the ballot theorem states that for the random walk $R$ with interchangeable and non-negative increments which starts at $0$ (i.e. $R(0)=0$), the probability that $R(m)<m$ for each $m=1,2,\ldots,n$, conditionally on $R(n)=k$, is equal to $\frac{k}{n}$.
More precisely, for us, the following result will play the key role. The result of this type was proved independently by \emph{Dwass} (for cyclically interchangeable random variables) and by \emph{Tak\'acs} (in less general case, for interchangeable random variables) in 1962. (appearing in the same issue of Annals of Mathematical Statistics), for details see \cite{Tak},\cite{Dwa} and for historical overview of the named results see \cite{AB}.
\begin{lemma}
Let $\xi_1,\ldots,\xi_n$ be the cyclically interchangeable random variables with values in the set $\{\ldots,-3,-2,-1,0,1\}$. Let $R(i)=\xi_1+\cdots+\xi_i$, $1\leq i\leq n$, $R(0)=0$. Then for each $0\leq k\leq n$
\begin{equation}\label{2.eq.1}
\P(R(i)>0~\textrm{for each  $1\leq i\leq n$}~|~R(n)=k)=\frac{k}{n}~.
\end{equation}
\end{lemma}
Let us notice that from Lemma 2.3. it follows that if we have a skip-free random walk and we know it's position at some instant $n$ and that position is some $k$, we are able to calculate the exact probability that this random walk stayed under the position $0$ (or above the position $0$, depending on which skip-free random walk we observe, the right or the left continuous one) and that probability is equal to $\frac{k}{n}$. It is also important for our problem to mention that the assumptions used in the above lemma cannot be removed - it is necessary that the variables take values in $\Z$ and that they are bounded from one side, i.e. that we can control the jumps on one side.\\
\\
To prove Lemma 2.3. we need the following result, again for details see \cite{Tak}.
\begin{lemma}
Let $\varphi(u)$, $u=0,1,2,\ldots$ be a nondecreasing function for which $\varphi(0)=0$ and
$\varphi(t+u)=\varphi(t)+\varphi(u)$, for $u=0,1,2,\ldots$,
where $t\in\N$. Define $$\delta(u)=\left\{%
\begin{array}{ll}
    1, & \hbox{$v-\varphi(v)>u-\varphi(u)$ za $v>u$;} \\
    0, & \hbox{otherwise.} \\
\end{array}%
\right.$$ Then
\begin{equation}\label{2:eq.3}\sum_{u=1}^{t}\delta(u)=\left\{%
\begin{array}{ll}
    t-\varphi(t), & \hbox{$0\leq\varphi(t)\leq t$;} \\
    0, & \hbox{$\varphi(t)\geq t$.} \\
\end{array}%
\right.\end{equation}
\end{lemma}
Let us now prove Lemma 2.3.\\
\\
We observe the random variables  $\gamma_1:=1-\xi_1$, $\gamma_2:=1-\xi_2$,$\ldots$ (instead of $\xi_1$,$\xi_2$,$\ldots$) and the random walk $R$  defined as $R(i)=\gamma_1+\ldots+\gamma_i$, $i\geq 1$, $R(0)=0$. It is a nondecreasing random walk and it's increments $\gamma_i$ are integer and cyclically interchangeable variables. We will show that
$$\P(R(u)\leq u~:~0\leq u\leq n|R(n))=\left\{%
\begin{array}{ll}
    1-\frac{R(n)}{n}, & \hbox{$0\leq R(n)\leq n$;} \\
    0, & \hbox{\textrm{otherwise.}} \\
\end{array}%
\right .$$
First we associate a new process $(R^{*}(u)~:~0\leq u<\infty)$ on $(0,\infty)$ to the process $(R(u)~:~0\leq u\leq n)$ such that
$R^{*}(u)=R(u)$, for $0\leq u\leq n$ and
$R^{*}(n+u)=R^{*}(n)+R^{*}(u)$, for $u\geq 0$. We define
$$\delta(u)=\left\{%
\begin{array}{ll}
    1, & \hbox{if $v-R^{*}(v)\geq u-R^{*}(u)$ for each $v\geq u$;} \\
    0, & \hbox{otherwise.} \\
\end{array}%
\right.$$
Then $\delta(u)$ is a random variable and has the same distribution for each $u\geq 0$.\\
\\
Now we have
 \begin{align*} &\P(R(u)\leq u~,~0\leq u\leq n|R(n))=\P(R^{*}(u)\leq u~,~u\geq 0|R(n))\\
&=\E[1_{\{v-R^{*}(v)\geq 0~,~\forall v\geq 0\}}|R(n)]=\E[\delta^{*}(0)|R(n)]=\\
&=\frac{1}{n}\cdot\sum_{u=1}^{n}\E[\delta^{*}(u)|R(n)]=\E[\frac{1}{n}\cdot\sum_{u=1}^{n}\delta^{*}(u)|R(n)]\\
&=\left\{%
\begin{array}{ll}
    1-\frac{R(n)}{n}, & \hbox{$0\leq R(n)\leq n$;} \\
    0, & \hbox{otherwise,} \\
\end{array}%
\right.
\end{align*}
using the fact that $\delta^{*}(u)$, $u\geq 0$, conditional on the position of $R(n)$, is equally distributed as $\delta(u)$, $0\leq u\leq n$, and using the result of Lemma 2.3. $\blacksquare$\\
\\
This is the proof which follows the approach used in \cite{Tak} and is also suitable for the continuous time risk process. But this type of result can also be proved following slightly different approach, in the same way the classic ballot theorems were proved - i.e. using some kind of a combinatorial formula. More precisely, we can use the following, for similar approach see \cite{Lam}.
\begin{lemma}(\textbf{combinatorial formula}) Let $(x_1,\ldots ,x_n)$ be a finite sequence of integers with values greater or equal to $-1$ and let
$\sum_{i=1}^{n} x_i=-k$. Let $\sigma_i(x)$ be a cyclic permutation of
$x$ which starts with $x_i$, i.e. $\sigma_i(x)=(x_i,x_{i+1},\ldots ,x_n,x_1,\ldots x_{i-1})$, $i\in\{1,2,\ldots ,n\}$.
Then there are exactly
$k$ different indices $i\in\{1,2,\ldots ,n\}$ such that
$$\sum_{l=1}^{j}(\sigma_i(x))_l> -k~,~~\textrm{for each}~j=1,\ldots ,n-1$$
and
$$\sum_{l=1}^{n}(\sigma_i(x))_l=-k~,$$
i.e. there are exactly $k$ different permutations $\sigma_i(x)$ of the sequence $x$ such that the first sum of the members of the sequence $\sigma_i(x)$ that is equal to $-k$ is the sum of all members of the sequence $\sigma_i(x)$.
\end{lemma}
\textit{Proof.}\\
\\
We observe partial sums $s_j=\sum_{i=1}^j x_i$, $1\leq j\leq n$, $s_0:=0$, and find the lowest one - let that be $s_m$ (i.e. $m$ is the lowest index such that $s_m=\min_{1\le j \le n} s_j$). Now we take the cyclic permutation $\sigma_m(x)$, i.e. the one that starts with $x_{m+1}$;
$\sigma_m(x)=(x_{m+1},x_{m+2},\ldots ,x_m)$. The overall sum of the sequence is $-k$, so $\sigma_m(x)$ hits $-k$ for the first time at the time instant $n$. For $j=1,2,\ldots ,k$ let $t_j$ be the first time of hitting the level $-j$, i.e. $t_j:=\min\{i\geq
1~:~s_i=-j\}$. Now again we can see that $\sigma_i(x)$ hits $-k$ for the first time at the time instant $n$ if and only if $i$ is one of the $t_j$-s, $j=1,2,\ldots ,k$, which proves our formula. $\blacksquare$\\
\\
Let us now observe the random walk $R(j)=\sum_{i=1}^{j} \xi(i)$, $R(0)=0$, with cyclically interchangeable increments $(\xi(1),\ldots ,\xi(n))$. Let $T_{-k}$ be the first time that $R$ reaches the level $-k$ and $T^{(i)}_{-k}$ the first time when the random walk with increments $\sigma_i(\xi)$ reaches $-k$ for the first time. Using Lemma 2.5., we have
\begin{align*}
n\cdot\P(T_{-k}=n|R(n)=-k)&=\sum_{i=1}^{n}
\E[1_{\{T_{-k}=n\}}|R(n)=-k]\\
&=\E\big(\sum_{i=1}^{n}1_{\{T_{-k}=n\}}|R(n)=-k\big)\\
&=\E\big(\sum_{i=1}^{n}1_{\{T^{(i)}_{-k}=n\}}|R(n)=-k\big)\\
&=k~,
\end{align*}
where in the second line from the end we used that the increments of the random walk $R$ are cyclically interchangeable and in the last line the combinatorial formula, i.e. the fact that there are exactly $k$ permutations of the increments of the random walk $R$ which hit the level $-k$ for the first time at the time instant $n$.
This is the \textit{Kemperman's formula} or \emph{the hitting time theorem}, and since we will use it only for independent random variables, i.e. the increments of the random walk, we will rephrase it in the less general form than the one we proved above.
\begin{lemma} Let $R$ be the upwards skip-free random walk starting at $0$ (i.e. $R(0)=0$) and $\tau(k)\in\{0,1,2,\ldots\}$ the first time that the random walk $R$ crosses the level $k>0$. Then
\begin{equation}\label{2:eq.4}n\cdot\P(\tau(k)=n)=k\cdot\P(R(n)=k)~,~n\geq 1~.
\end{equation}
\end{lemma}

\section{Main results for the discrete time ruin process} \label{sec-3}
Let $C^1$and $C^2$ be the independent and nondecreasing random walks, i.e.
\begin{equation}\label{3:eq.1}
C^i(n)=U^i(1)+\cdots +U^i(n)~,~~n\geq 1~,
\end{equation}
for
\begin{equation}\label{3:eq.2}
U^i(j)\sim \left(%
\begin{array}{ccccc}
  0 & 1 & 2 & 3 & \ldots \\
  p_0 & p_1 & p_2 & p_3 & \ldots
\end{array}%
\right)
\end{equation}
for some $p_k\geq 0$, $k\geq 0$, $\sum_{k=0}^\infty p_k=1$ and $j\in\{1,2,\ldots\}$, $i=1,2$. We define
\begin{equation}\label{3:eq.3}
C=C^1+C^2~.
\end{equation} Let us assume that $\E C^i(1)<\infty$ or, equivalently, $\E U^i(1)<\infty$, $i=1,2$.\\
\\
We define the \emph{the discrete time risk process with the unit drift} by
\begin{equation}\label{3:eq.3}
X(n)=n-C(n)~,~~n\geq 0~,
\end{equation}
which means that $X$ is the upwards skip-free random walk. Let us further assume that, using the notation  $\mu:=\E C(1)$, $\mu_i:=\E C^i(1)$, $i=1,2$,
\begin{equation}\label{3:eq.4}
\E X(1)=1-\E C(1)=1-\mu=1-(\mu_1+\mu_2)>0~.
\end{equation}
For $\widehat{X}=-X$, we also define
\begin{equation*}\widehat{S}(n):=\max_{0\leq s\leq n} \widehat{X}(s)~,
\end{equation*}
\begin{equation*}
\widehat{S}(\infty):=\max_{s\geq 0} \widehat{X}(s)
\end{equation*}
and
\begin{equation*}\widehat{\tau_x}=\inf\{n\geq 0~:~\widehat{X}(n)>x\}~,
\end{equation*}
the first time that the dual random walk $\widehat{X}$ crosses the level  $x\in\N$,
\begin{equation*}
\Delta C^i(n)=C^i(n)-C^i(n-1)~~\textrm {for} ~~i=1,2
\end{equation*}
and
\begin{equation*}
\Delta C(n)=C(n)-C(n-1)~,
\end{equation*}
$n\geq 1$, the jumps of the random walks $C^1$, $C^2$ and $C$.
Using the linearity of the expectation, the fact that the increments of the random walk are independent and equally distributed and the standard induction procedure, we can see that the following result is valid.
\begin{lemma}
\begin{equation}\label{3:eq.5}\E\big( \sum_{n=0}^{\infty} \mathcal{H}(n,\omega,\Delta
C_n^i(\omega))\big) =\E\big(\int_{(0,\infty)}
\mathcal{H}(n,\omega,\varepsilon)dF_i(\varepsilon)\big)~,
\end{equation}
where $\mathcal{H}$ is a non-negative function and $F_i$ the distribution function of the increments of the random walk $C^i$, $i=1,2$.
\end{lemma}
For $y\leq 0$ and $x>0$ we define
$$\mathcal{H}(n,\omega,\varepsilon_i)=1_{\{\widehat{X}(n-1)=y,\widehat{S}(n-1)\leq
0\}}\cdot 1_{\{\varepsilon_i=x+1-y\}}~.$$
Then, using Lemma 3.1., it follows
\begin{align*} &\E\big( \sum_{n=0}^{\infty} \mathcal{H}(n,\omega,\Delta
C^i(n)(\omega))\big)\\
&=\sum_{n=1}^{\infty} \P(\widehat{X}(n-1)=y,\widehat{S}(n-1)\leq
0,\Delta C^i(n)=x+1-y)\\
&=\P(\widehat{\tau_0}<\infty,\widehat{X}(\widehat{\tau_0}-1)=y,
\widehat{X}(\widehat{\tau_0})\geq x,\Delta
C^i(\widehat{\tau_0})=x+1-y)~.
\end{align*}
Let us mention that the inequality $\widehat{X}(\widehat{\tau_0})\geq x$ appearing in the last line is the result of the fact that in the discrete time case the components of the random walk $C$ may jump simultaneously (unlike in the continuous time case when modelling the risk process with the spectrally negative L\'evy process). Since the components of the random walk $C$ are nondecreasing, they can only increase the supremum of the overall risk process and the drift decreases it for a unit at each time instant - so we have $\Delta
C^i(\widehat{\tau_0})=(x-y)+1$ in the last line.\\
\\
On the other side, using Lemma 2.3. and Lemma 2.6., we have\\
\\
$\E\big(\sum_{n=1}^{\infty}\int_{(0,\infty)}
\mathcal{H}(n,\omega)dF_i(\varepsilon_i)\big)\\
\\=\sum_{n=1}^{\infty}\E\big(\int_{(0,\infty)}
1_{\{\widehat{X}(n-1)=y,\widehat{S}(n-1)\leq 0\}}\cdot
1_{\{\varepsilon_i=x+1-y\}}dF_i(\varepsilon_i)\big)\\
\\=\sum_{n=1}^{\infty}\P(\widehat{X}(n-1)=y,\widehat{S}(n-1)\leq 0)\cdot\P(C^i(1)=x+1-y)\\
\\=\sum_{n=1}^{\infty}\P(\widehat{S}(n-1)\leq 0|\widehat{X}(n-1)=y)\cdot\P(\widehat{X}(n-1)=y)\cdot\P(C^i(1)=x+1-y)\\
\\=\P(C^i(1)=x+1-y)\cdot\sum_{n=1}^{\infty} \P(\max_{0\leq m\leq
n-1} \widehat{X}(m)\leq
0|\widehat{X}(n-1)=y)\cdot\P(\widehat{X}(n-1)=y)\\
\\=\P(C^i(1)=x+1-y)\cdot\sum_{n=1}^{\infty} \P(\widehat{X}(m)\leq
0~,~ \forall 0\leq m\leq n-1|\widehat{X}(n-1)=y)\cdot\P(\widehat{X}(n-1)=y)\\
\\=\P(C^i(1)=x+1-y)\cdot\sum_{n=1}^{\infty}\P(C(m)\leq m~:~\forall 0\leq m\leq n-1|\widehat{X}(n-1)=y)\cdot\P(\widehat{X}(n-1)=y)\\
\\=\P(C^i(1)=x+1-y)\cdot\sum_{n=1}^{\infty}\P(C(m)\leq m~:~ \forall 0\leq m\leq
n-1|C(n-1)=y+(n-1))\cdot\P(\widehat{X}(n-1)=y)\\
\\=\P(C^i(1)=x+1-y)\cdot\sum_{n=1}^{\infty} (1-\frac{y+(n-1)}{n-1})\cdot\P(\widehat{X}(n-1)=y)\\
\\=\P(C^i(1)=x+1-y)\cdot\sum_{n=1}^{\infty} \frac{1}{n-1}\cdot(-y)\cdot\P(\widehat{X}(n-1)=y)\\
\\=\P(C^i(1)=x+1-y)\cdot\sum_{n=1}^{\infty} \frac{1}{n-1}\cdot(n-1)\cdot\P(\widehat{\tau}(y)=n-1)\\
\\=\P(C^i(1)=x+1-y)~.$\\
\\
Let us notice that in the last line we again used the fact that $\widehat{X}$ is a downwards skip-free random walk which can only take unit steps to go downwards, i.e. it has to hit each level $y=-k\leq 0$ it crosses, so
$\P(\widehat{\tau}(y)<\infty)=1$ for $y\leq 0$.\\
\\
Let us further notice that the above result can be generalized for finitely many random walks $C^1,C^2,\ldots ,C^m$, for some $m\in\N$, in the same way. So we have the following result.
\begin{thm} Let $C^1, C^2,\ldots ,C^m$, $m\in\N$, be independent random walks with nondecreasing increments (defined as in (\ref{3:eq.1}) and (\ref{3:eq.2})) and
\begin{equation}\label{3:eq.7}
X(n)=n-(C^1+\cdots
+C^m)(n)~,~n\geq 0~.
\end{equation}
Let $\widehat{X}(0)=0$ for the dual of the random walk $X$  (i.e. $\widehat{X}=-X$) and let us assume that $1>\mu=\mu^1+\cdots +\mu^m$, for $\mu:=\E C(1)$ and $\mu^i:=\E C^i(1)$, $i=1,2,\ldots,m$. Then for $y\leq 0$ and $x>0$
\begin{equation}\label{3:eq.8}\P(\widehat{\tau_0}<\infty,\widehat{X}(\widehat{\tau_0}-1)=y,
\widehat{X}(\widehat{\tau_0})\geq x,\Delta
C^i(\widehat{\tau_0})=x+1-y)=\P(C^i(1)=x+1-y)~.
\end{equation}
\end{thm}
Summing for all $y\leq 0$, we derive the following result.
\begin{cor} For the random walks and the assumptions defined as in the Theorem 3.2.,
\begin{equation}\label{3:eq.9}\P(\widehat{\tau_0}<\infty, \widehat{X}(\widehat{\tau_0})\geq
x,\Delta C^i(\widehat{\tau_0})\geq x+1)=\P(C^i(1)\geq x+1)~.
\end{equation}
\end{cor}
Let us now look at \emph{the discrete risk model with the perturbation}, i.e. the random walk $X$ such that
\begin{equation}\label{3:eq.10}
X(n)=-C(n)+Z(n)~,~ n\geq 1~,\end{equation}
where $C$ is the random walk with nondecreasing increments and $Z$ the upwards skip-free random walk. In other words, we have
\begin{equation}\label{3:eq.11}
Z(1)=\xi_Z(i)\sim\left(%
\begin{array}{ccccc}
  \ldots & -2 & -1 & 0 & 1 \\
  \ldots & q_2 & q_1 & q_0 & \rho \\
\end{array}%
\right)
\end{equation}
and
\begin{equation}\label{3:eq.12}
C(1)=\xi_C(i)\sim\left(%
\begin{array}{ccccc}
   0 & 1 & 2 & 3 & \ldots \\
   p_0& p_1 & p_2 & p_3 & \ldots \\
\end{array}%
\right)~~,
\end{equation}
for some $\rho,~q_j,~p_j\geq 0$ such that $\sum_{j=0}^\infty q_j=\sum_{j=0}^\infty p_j=1$.
We assume that
\begin{equation}\label{3:eq.13}
\E(X(1))>0~, ~~\textrm{i.e}~~\E C(1)<\E Z(1)~.
\end{equation}
$X$ is obviously the upwards skip-free random walk and the dual process, $\widehat{X}=-X$, is the downwards skip-free random walk such that $\widehat{X}\to -\infty$. Furthermore, we can rewrite $\widehat{X}$ so that
$$\widehat{X}(n)=\sum_{i=1}^{n}\xi_{\widehat{X}}(i)=\sum_{i=1}^{n}(\xi_W(i)-1)=\sum_{i=1}^n
\xi_W(i)-n=:W(n)-n~,~~n\geq 0~,$$
for $\xi_W(i):=\xi_{\hat{X}}(i)+1$, so we have that $\P(W(1)=k)=\P(\widehat{X}(1)=k-1)$, $k\geq 0$.\\
\\
Since in the discrete time case the random walks $C$ and $Z$ jump at the same time, if $\widehat{X}$ was in some position $y\leq 0$
at the time instant just before it crossed the level $0$ and in the position $x>0$ when it crossed the level $0$, the event  \{$C$ caused the jump of the process $\widehat{X}$ over the level $0$\} can be written as
$$\{\Delta C(\widehat{\tau}_0)\geq -y+1+x, \Delta Z(\widehat{\tau}_0) =-1\}\cup\{\Delta C(\widehat{\tau}_0)\geq -y+x, \Delta Z(\widehat{\tau}_0) \geq 0\}~.$$
So, for $y\leq 0$ and $x>0$ we have
\begin{align*} &\P(\widehat{\tau}_0<\infty,\widehat{X}(\widehat
{\tau}_0-1)=y,\widehat{X}(\widehat{\tau}_0)\geq
x,\textrm{$C$ caused the jump of the process $\widehat{X}$ over the level $0$ })\\
&=\P(\widehat{\tau}_0<\infty,\widehat{X}(\widehat{\tau}_0-1)=y,\widehat{X}(\widehat{\tau}_0)\geq
x,\Delta C(\widehat{\tau}_0)\geq -y+1+x,\Delta
Z(\widehat{\tau}_0)=-1)\\
&+\P(\widehat{\tau}_0<\infty,\widehat{X}(\widehat{\tau}_0-1)=y,\widehat{X}(\widehat{\tau}_0)\geq
x,\Delta C(\widehat{\tau}_0)\geq -y+x,\Delta
Z(\widehat{\tau}_0)\geq 0)~.
\end{align*}
Let us define
$$\mathcal{H}(n,\omega,\varepsilon,\eta)=1_{\{\widehat{X}(n-1)=y,\widehat{S}(n-1)\leq
0\}}\cdot 1_{\{\varepsilon=x+1-y\}}\cdot 1_{\{\eta=-1\}}~.$$
Using the Lemma 3.1., Lemma 2.3. and Lemma 2.6. we have
\begin{align*}
&\P(\widehat{\tau}_0<\infty,\widehat{X}(\widehat{\tau}_0-1)=y,\widehat{X}(\widehat{\tau}_0)\geq
x,\Delta C(\widehat{\tau}_0)\geq -y+1+x,\Delta
Z(\widehat{\tau}_0)=-1)\\
&=\sum_{n=1}^{\infty}\P(\widehat{X}(n-1)=y,\widehat{S}(n-1)\leq
0)\cdot\P(C(1)\geq x+1-y)\cdot\P(Z(1)=-1)\\
&=\P(C(1)\geq x+1-y)\cdot\P(Z(1)=-1)\cdot\sum_{n=1}^{\infty}\P(\widehat{S}(n-1)\leq
0|\widehat{X}(n-1)=y)\cdot\P(\widehat{X}(n-1)=y)\\
&=\P(C(1)\geq x+1-y)\cdot\P(Z(1)=-1)\cdot\sum_{n=1}^{\infty}\P(W(t)\leq t~,~0\leq t\leq
n-1|\widehat{X}(n-1)=y)\\
& \cdot\P(\widehat{X}(n-1)=y)\\
&=\P(C(1)\geq x+1-y)\cdot\P(Z(1)=-1)\cdot\sum_{n=1}^{\infty}(1-\frac{y+(n-1)}{n-1})\cdot\P(\widehat{X}(n-1)=y)=\\
&=\P(C(1)\geq x+1-y)\cdot\P(Z(1)=-1)\cdot\sum_{n=1}^{\infty}\frac{1}{n-1}\cdot(-y)\cdot\P(\widehat{X}(n-1)=y)=\\
&=\P(C(1)\geq
x+1-y)\cdot\P(Z(1)=-1)\cdot\sum_{n=1}^{\infty}\P(\widehat{\tau}_{y}=n-1)=\\
&=\P(C(1)\geq x+1-y)\cdot\P(Z(1)=-1)~.\end{align*}
It is obvious that in the result above we used the crucial argument that $\widehat{X}$ is the downwards skip-free random walk, which means that it can only use unit steps to go downwards, i.e. it visits each level $y\leq 0$.\\
\\
We can use the same calculation on the second addend,  $\P(\widehat{\tau}_0<\infty,\widehat{X}(\widehat{\tau}_0-1)=y,\widehat{X}(\widehat{\tau}_0)\geq
x,\Delta C(\widehat{\tau}_0)\geq -y+x,\Delta
Z(\widehat{\tau}_0)\geq 0)$, so we have the following result.
\begin{thm} Let $C$ and $Z$ be the random walks defined as in (\ref{3:eq.11}) and (\ref{3:eq.12}) and $X$ the discrete time perturbed risk process
\begin{equation}\label{3:eq.14}
X(n)=-C(n)+Z(n)~,~ n\geq 1~,\end{equation} and let us assume that $\E(X(1))<0$, i.e. $\E C(1)<\E Z(1)$. Then
\begin{align*}\P(\widehat{\tau}_0<\infty,&\widehat{X}(\widehat
{\tau}_0-1)=y,\widehat{X}(\widehat{\tau}_0)\geq x,\textrm{$C$ caused the jump of the process $\widehat{X}$ over the level $0$})\\
&=\P(C(1)\geq
x+1-y)\cdot\P(Z(1)=-1)+\P(C(1)\geq x-y)\cdot\P(Z(1)\geq
0)~.\end{align*}\end{thm}
Let us notice that the similar result can be derived if we define $C$ as the sum of the $m\in\N$ independent nondecreasing random walks $C^1,\ldots,C^m$, $C=C^1+\cdots +C^m$.\\
\\
Let us further notice that the "problem" of the simultaneous jumps of the random walks $C$ and $Z$, which is characteristic for the discrete time processes and differs from the continuous time version of the same problem, may be overcome if we observe a natural connection between these two types of models - i.e. the compound Poisson processes. For some generalizations in this way and similar results connected to ruin probability see \cite{GTV}.\\
\\
\textbf{Acknowledgement:} This work has been supported in part by Croatian Science Foundation under
the project 3526.

\end{document}